\documentclass[12 pt,twosided,reqno]{amsart}
\usepackage{amscd,amssymb}
\usepackage[arrow,matrix]{xy}
\usepackage{graphicx}
\usepackage{epsfig}
\theoremstyle{plain}
\newtheorem{theorem}[subsection]{Theorem}

\newtheorem{proposition}[subsection]{Proposition}
\newtheorem{cor}[subsection]{Corollary}

\theoremstyle{definition}
\newtheorem{remark}[subsection]{Remark}
\newtheorem{definition}[subsection]{Definition}
\newtheorem{example}[subsection]{Example}

\numberwithin{equation}{section} 

\def\deg{\hbox{deg}}
\def\br{\hbox{br}}

\def\Div{\hbox{Div}}

\begin{document}
\title [Recurrence Relation For HOMFLY polynomial\\ And Rational Specializations]
{Recurrence Relation For HOMFLY polynomial\\ And rational Specializations }
\thanks{\emph{keywords and phrases:} HOMFLY polynomial, Alexander-Conway polynomial, simple braids, Fibonacci recurrence}
\thanks{This research is partially supported by Higher Education Commission, Pakistan.\\
2010 AMS classification: Primary 57M27, 57M25; Secondary 20F36.}
\author{  REHANA ASHRAF$^{1}$,\,\,BARBU BERCEANU$^{1,2}$}

\address{$^{1}$Abdus Salam School of Mathematical Sciences,
 GC University, Lahore-Pakistan.}
\email {rashraf@sms.edu.pk}
\address{$^{2}$Institute of Mathematics Simion Stoilow, Bucharest-Romania (permanent address).}
\email {Barbu.Berceanu@imar.ro}
  \maketitle
 \pagestyle{myheadings} \markboth{\centerline {\scriptsize
 REHANA ASHRAF,\,\,\,BARBU BERCEANU  }} {\centerline {\scriptsize
Recurrence Relation For HOMFLY Polynomial And Rational Specializations}}
\begin{abstract}Turning the skein relation for HOMFLY into a Fibonacci recurrence,
we prove that there are only three rational specializations of HOMFLY polynomial: Alexander-Conway, Jones,  and a new one. Using the recurrence relation, we find  general and relative expansion formulae and rational generating functions for Alexander-Conway polynomial and the new polynomial, which reduce the computations to closure of  simple braids, a subset of square free braids;  HOMFLY polynomials of these simple braids are also computed. Algebraic independence of these three polynomials is proved.
\end{abstract}
\section{\textbf{Introduction}} \label{sec1}
HOMFLY polynomial (see \cite{HOM:FLY}, \cite{Lickorish:97}) is a function $P:\{\hbox{oriented link diagrams}\}\longrightarrow\mathbb{C}[l^{\pm 1},m^{\pm 1}]$ satisfying:\\
\indent \textbf{\emph{a)}} \textbf{Normalization:} $P(\bigcirc)=1;$\\
\indent \textbf{\emph{b)}} \textbf{Skein Relation:} whenever three oriented link diagrams
$L_{+} ,L_{-}$ and  $  L_{0}$ are the same,  except in the
neighborhood of one crossing, where

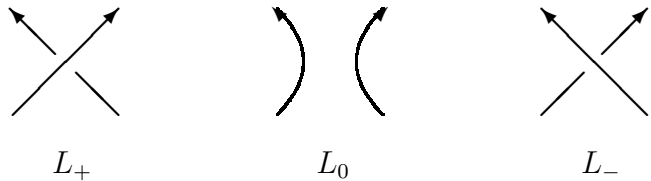
\begin{center}
\begin{picture}(360,100)
\thicklines
 \put(75,18){$L_{+}$}
  \put(175,18){$L_{0}$}
   \put(275,18){$L_{-}$}

 \qbezier(160,40)(180,60)(160,80)
  \qbezier(200,40)(180,60)(200,80)
  \put(162,77){\vector(-1,1){4}}
  \put(198,77){\vector(1,1){4}}
  \put(60,40){\vector(1,1){40}}
    \put(283,63){\vector(1,1){17}}
   \put(300,40){\vector(-1,1){40}}
     \put(76,63){\vector(-1,1){17}}
     \put(100,40){\line(-1,1){16}}
        \put(260,40){\line(1,1){16}}
\end{picture}
\end{center}
then\\
  $$lP(L_{+})+l^{-1}P(L_{-})+mP(L_{0})=0.$$
HOMFLY polynomial is in fact a link invariant which generalizes both Alexander-Conway and Jones polynomials. The specializations:\\
\indent 1)  $l=\iota t^{-1},m=\iota(t^{-\frac{1}{2}}-t^{\frac{1}{2}})$ gives
the Jones polynomial $V(L)$ and  skein relation
$$tV(L_{-})-(t^{-\frac{1}{2}}-t^{\frac{1}{2}})V(L_{0})-t^{-1}V(L_{+})=0
,\, \,\mbox{\rm {respectively}}$$
\indent 2)  $l=\iota ,m=\iota(t^{\frac{1}{2}}-t^{-\frac{1}{2}})$ gives the
Alexander-Conway polynomial $\nabla(L)$ and  skein relation
$$\nabla(L_{-})-(t^{\frac{1}{2}}-t^{-\frac{1}{2}})\nabla(L_{0})-\nabla(L_{+})=0.$$

\medskip

In this paper we transform the skein relation into a Fibonacci type
recurrence for HOMFLY polynomial of  closed  braids. The generators $x_{i}\,(i=1,\ldots, n-1)$ of the braid group
$\mathcal{B}_{n}$ are classical Artin generators
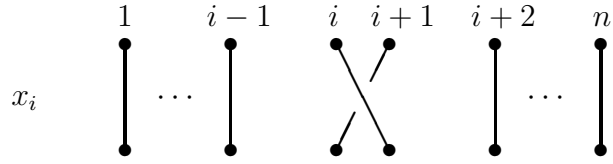
\begin{center}
\begin{picture}(360,80)
\thicklines
 \put(90,60){$1$}              \put(90,10){$\bullet$}    \put(93,13){\line(0,1){40}}
  \put(90,50){$\bullet$}       \put(270,60){$n$}         \put(270,10){$\bullet$}
\put(273,13){\line(0,1){40}}   \put(270,50){$\bullet$}   \put(124,60){${i-1}$}
\put(130,10){$\bullet$}     \put(133,13){\line(0,1){40}} \put(130,50){$\bullet$}
\put(224,60){$i+2$}         \put(230,10){$\bullet$}      \put(233,13){\line(0,1){40}}
\put(230,50){$\bullet$}     \put(170,60){$i$}            \put(170,10){$\bullet$}
\put(170,50){$\bullet$}     \put(186,60){$i+1$}          \put(190,10){$\bullet$}
\put(190,50){$\bullet$}    \put(193,13){\line(-1,2){20}} \put(173,13){\line(1,2){7}}
\put(193,53){\line(-1,-2){7}} \put(105,30){$\cdots$}     \put(245,30){$\cdots$}
 \put(50,30){$x_{i}$} 
\end{picture}
\end{center}
The braid group $\mathcal{B}_{n}$ is fixed, and also the sequence of
generators $(x_{i_{1}},x_{i_{2}},\ldots, x_{i_{k}})$. Let us denote
by $P_{n}(a_{1},\ldots, a_{k})$ HOMFLY polynomial of the closure of
the $n$-braid $\beta=x_{i_{1}}^{a_{1}}x_{i_{2}}^{a_{2}}\ldots
x_{i_{k}}^{a_{k}}$. The first result of the paper is
\begin{theorem}
\label{recurrence: theorem}
For any $a_{1},\ldots,a_{k}$ $\in \mathbb{Z}$ and any $j=1,\ldots,k $ the recurrence holds\\
$$P_{n}(a_{1},\ldots,a_{j}+2,\ldots,a_{k})+mlP_{n}(a_{1},\ldots,a_{j}+1,\ldots,a_{k})+l^{2}P_{n}(a_{1},\ldots,a_{j},\ldots,a_{k})=0.$$
\end{theorem}

\medskip

Using the terminology of \cite{nizami:08}, the sequence
$\{P_{n}(a_{1},\ldots, a_{k})\}_{a_{*}\in\mathbb{Z}^{k}}$ is \emph{a
multiple Fibonacci sequence} with parameters $(-ml,-l^{2})$. Make
the change of variable $t\rightarrow s^{-2}$, Alexander-Conway and
Jones specializations give two recurrences;  in the next formulae
$\nabla_{n}(a_{1},\ldots,a_{k})$, $V_{n}(a_{1},\ldots,a_{k})$ (and also $D_{n}(a_{1},\ldots,a_{k})$)
stand for corresponding polynomials of the closure of $n$-braid
$\beta=x_{i_{1}}^{a_{1}}x_{i_{2}}^{a_{2}}\ldots x_{i_{k}}^{a_{k}}.$

\begin{cor}\label{cor:1}
Alexander-Conway polynomial satisfies the multiple Fibonacci recurrence:\\
$$\nabla_{n}(a_{1},\ldots,a_{j}+2,\ldots,a_{k})=(s^{-1}-s)\nabla_{n}(a_{1},\ldots,a_{j}+1,\ldots,a_{k})+\nabla_{n}(a_{1},\ldots,a_{j},\ldots,a_{k}).$$
\end{cor}
\medskip
\begin{cor}\mbox{\rm {(see} \cite{Barbu-Nizami:09})}\,Jones polynomial satisfies the multiple Fibonacci recurrence:
$$V_{n}(a_{1},\ldots,a_{j}+2,\ldots,a_{k})=(s^{3}-s)V_{n}(a_{1},\ldots,a_{j}+1,\ldots,a_{k})+s^{4}V_{n}(a_{1},\ldots,a_{j},\ldots,a_{k}).$$
\end{cor}
A "relative" example of multiple Fibonacci sequence is given by HOMFLY polynomials of the link diagrams $\mathcal{L}(R_{k_{1}},\ldots,R_{k_{p}})$, where  there are $p$ subdiagrams  and the diagram with index $j$ contains the braid $x_{1}^{k_{j}}$, ($k_{1},\ldots,k_{p}$ are arbitrary integers)

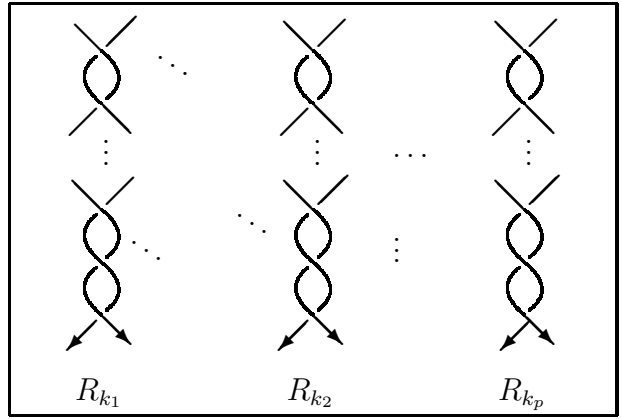
\begin{center}
\begin{picture}(360,160)
\thicklines
 \qbezier(80,140)(92,130)(82,122)    \qbezier(78,138)(68,130)(80,120)
 \put(80,140){\line(-1,1){10}}       \put(82,142){\line(1,1){11}}
 \put(80,120){\line(1,-1){11}}       \put(78,118){\line(-1,-1){10}}
 \qbezier(80,80)(92,70)(82,62)       \qbezier(78,78)(68,70)(80,60)
 \qbezier(80,60)(92,50)(82,42)       \qbezier(78,58)(68,50)(80,40)
 \put(80,80){\line(-1,1){11}}        \put(82,82){\line(1,1){10}}
 \put(80,40){\vector(1,-1){11}}      \put(78,38){\vector(-1,-1){11}}
 \qbezier(160,140)(172,130)(162,122) \qbezier(158,138)(148,130)(160,120)
 \put(160,140){\line(-1,1){11}}      \put(160,120){\line(1,-1){11}}
 \put(162,142){\line(1,1){10}}       \put(158,118){\line(-1,-1){10}}
 \qbezier(160,80)(172,70)(162,62)    \qbezier(158,78)(148,70)(160,60)
 \qbezier(160,60)(172,50)(162,42)    \qbezier(158,58)(148,50)(160,40)
 \put(160,80){\line(-1,1){11}}       \put(162,82){\line(1,1){11}}
 \put(160,40){\vector(1,-1){11}}     \put(158,38){\vector(-1,-1){11}}
 \qbezier(240,140)(252,130)(242,122) \qbezier(238,138)(228,130)(240,120)
 \put(240,140){\line(-1,1){11}}      \put(240,120){\line(1,-1){11}}
 \put(242,142){\line(1,1){10}}       \put(238,118){\line(-1,-1){10}}
 \qbezier(240,80)(252,70)(242,62)    \qbezier(238,78)(228,70)(240,60)
 \qbezier(240,60)(252,50)(242,42)    \qbezier(238,58)(228,50)(240,40)
 \put(240,80){\line(-1,1){11}}       \put(242,82){\line(1,1){11}}
 \put(240,40){\vector(1,-1){11}}     \put(242,38){\vector(-1,-1){11}}
 \put(45,2){\line(1,0){230}}         \put(45,158){\line(1,0){230}}
 \put(45,2){\line(0,1){156}}         \put(275,2){\line(0,1){156}}
 \put(70,8){$R_{k_{1}}$}\put(150,8){$R_{k_{2}}$}\put(230,8){$R_{k_{p}}$}
 \put(80,97){\vdots}                 \put(160,97){\vdots} 
 \put(240,97){\vdots}                \put(190,100){\ldots}
 \put(190,60){\vdots}                \put(130,70){$\ddots$}
 \put(90,60){$\ddots$}               \put(100,130){$\ddots$}
 \put(2000,60){\vdots}               \end{picture}
 \end{center}
 and the rest of the diagram is fixed. HOMFLY polynomials of this family, denoted by $P_{\mathcal{L}}(k_{1},\ldots,k_{p})$, satisfy the same multiple Finonacci recurrence:
\begin{proposition} For any $k_{1},\ldots,k_{p}$ and $j\in\{1,2,\ldots,p\}$
$$P_{\mathcal{L}}(k_{1},\ldots,k_{j}+2,\ldots,k_{p})=-mlP_{\mathcal{L}}(k_{1},
\ldots,k_{j}+1,\ldots,k_{p})-l^{2}P_{\mathcal{L}}(k_{1},\ldots,k_{j},\ldots,k_{p}).$$
\end{proposition}

The roots of characteristic equations corresponding to these recurrences are $r_{1}=s^{-1},\,r_{2}=-s$ for Alexander-Conway and $r'_{1}=s^{3},\,r'_{2}=-s$ for Jones \cite{Barbu-Nizami:09}. It turns out that, up to constants, Alexander-Conway and Jones specializations are unique with the property that the roots of the characteristic equation are $s$ and $s^{m} ,m\in\mathbb{Z},\,m\neq1\,$; the case of double roots $r''_{1}=r''_{2}=s$ gives a new polynomial (see Theorem \ref{jones:conway}).\\
\indent\,\, Section 2 of the paper contains the proof of Theorem \ref{recurrence: theorem} and  precise definitions of rational and  normalized specializations (see Definition \ref{rat:ional} and Definition \ref{nor:malized}).\\
\indent\,\,In the next section two direct consequences of the recurrence relation for Alexander-Conway polynomial are given. The module of multiple Fibonacci sequences is a tensor power of a free module (see \cite{nizami:08}), therefore one can describe such sequences using a canonical basis: let us introduce the  Laurent polynomial\\
 $$C^{a}(s)=(-1)^{a}s^{a-1}+s^{1-a}.$$

\medskip

 \begin{theorem}
 \textbf{(Expansion formula)}
 \label{expansion:formula}
 The Alexander-Conway polynomial  of the closure of the $n$-braid $\beta=x_{i_{1}}^{a_{1}}x_{i_{2}}^{a_{2}}\ldots x_{i_{k}}^{a_{k}}$ is given by\\
 $$\nabla_{n}(x_{i_{1}}^{a_{1}} \ldots
x_{i_{k}}^{a_{k}})(s)=\Big(\frac{s}{s^2+1}\Big)^{k}\sum_{0\leq
j_{1},\ldots,j_{k}\leq1}C^{a_{1}+j_{1}}(s)\ldots
C^{a_{k}+j_{k}}(s)\nabla_{n}( x_{i_{1}}^{j_{1}} \ldots
x_{i_{k}}^{j_{k}})(s).$$
\end{theorem}
 Practitioners of  Burau representation of braid groups could compare the complexity of computation of Burau matrix and its determinant with the expansion formula. People working with Pretzel links and $2$-bridge knots  could find the next result familiar:
\begin{cor}\textbf{(Relative expansion formula)} The Alexander-Conway polynomial of the link $\mathcal{L}(R_{k_{1}},\ldots,R_{k_{p}})$ is $$\nabla_{\mathcal{L}}(k_{1},\ldots,k_{p})=\Big(\frac{s}{s^2+1}\Big)^{p}\sum_{0\leq
j_{1},\ldots,j_{p}\leq1}C^{k_{1}+j_{1}}(s)\ldots
C^{k_{p}+j_{p}}(s)\nabla_{\mathcal{L}}( {j_{1}}, \ldots,
{j_{p}})(s).$$
\end{cor}
\noindent For the next result $n$ and the sequence
$(x_{i_{1}},\ldots,x_{i_{k}})$ of generators of $\mathcal{B}_{n}$
are fixed. In order to compute the \emph{generating function} of Alexander-Conway polynomials of type $I_{*}=(i_{1},\ldots,i_{k})$, the formal series  in $t_{1},\ldots, t_{k}$:
$$\mathcal{AC}_{n}(t_{1},\ldots,t_{k})=\sum_{(a_{1},a_{2},\ldots,a_{k})\in\mathbb{Z}^{k}}\nabla_{n}(x_{i_{1}}^{a_{1}} \ldots x_{i_{k}}^{a_{k}})(s)t_{1}^{a_{1}}\ldots t_{k}^{a_{k}},$$
we introduce new polynomials in the formal variable $\tau:$\\

$q(\tau)=(1+s\tau)(1-s^{-1}\tau),\,Q_{0}(\tau)=1-(s^{-1}-s)\tau$ and
$Q_{1}(\tau)=\tau.$

\medskip

\begin{theorem}\textbf{(Generating function)}
\label{generating:function}
The generating function $\mathcal{AC}_{n}(t_{1},\ldots,t_{k})$ for Alexander-Conway polynomials is a rational function in $t_{1}.\ldots,t_{k}$:\\
$$\mathcal{AC}_{n}(t_{1},\ldots,t_{k})=q(t_{1})^{-1}\ldots q(t_{k})^{-1}\sum_{0\leq
j_{1},\ldots,j_{k}\leq1}Q_{j_{1}}(t_{1})\ldots
Q_{j_{k}}(t_{k})\nabla_{n}( x_{i_{1}}^{j_{1}}
\ldots x_{i_{k}}^{j_{k}}).$$
\end{theorem}
\noindent This section  contains also some computations of Alexander-Conway polynomial.

In  section 4 the degenerate case  ($r_{1}=r_{2}=s$) is analyzed. We find a new polynomial $D$ satisfying:\\
\indent \textbf{\emph{a)}}  \textbf{Normalization:} $D(\bigcirc)=1;$\\
\indent \textbf{\emph{b)}} \textbf{Skein Relation}:  $sD(L_{+})+s^{-1}D(L_{-})-2D(L_{0}) = 0$,\\
and corresponding expansion formulae
\begin{theorem}
\label{expansion:D}(\textbf{Expansion formula}) The D polynomial of the closure of the  $n$-braid $\beta=x_{i_{1}}^{a_{1}}x_{i_{2}}^{a_{2}}\ldots x_{i_{k}}^{a_{k}}$ is
$$D_{n}(x_{i_{1}}^{a_{1}}x_{i_{2}}^{a_{2}}\ldots x_{i_{k}}^{a_{k}})(s)=\sum_{0\leq
j_{1},\ldots,j_{k}\leq1}R^{[a_{1}]}_{j_{1}}(s)\ldots
R^{[a_{k}]}_{j_{k}}(s)D_{n}( x_{i_{1}}^{j_{1}} \ldots
x_{i_{k}}^{j_{k}})(s).$$
\end{theorem}
\medskip
\noindent (in the previous and next formulae: $R^{[a]}_{0}(s)=(1-a)s^{a},\,R^{[a]}_{1}(s)=as^{a-1}$),
\begin{cor}\textbf{(Relative expansion formula)} The $D$ polynomial of the link $\mathcal{L}(R_{k_{1}},\ldots,R_{k_{p}})$ is
$$D_{\mathcal{L}}(k_{1},\ldots,k_{p})=\sum_{0\leq
j_{1},\ldots,j_{p}\leq1}R^{[k_{1}]}_{j_{1}}(s)\ldots
R^{[k_{p}]}_{j_{p}} D_{\mathcal{L}}( {j_{1}}, \ldots
,{j_{p}})(s).$$

\end{cor}

\noindent and a rational \emph{generating function}
\begin{theorem}\label{generating:D}(\textbf{Generating function}) The generating function $\mathcal{D}_{n}(t_{1},\ldots ,t_{k})$ for the degenerate polynomials $D_{n}(x_{i_{1}}^{a_{1}}x_{i_{2}}^{a_{2}}\ldots x_{i_{k}}^{a_{k}})$ is a rational function in $t_{1},\ldots ,t_{k}$\\
$$\mathcal{D}_{n}(t_{1},\ldots,t_{k})=e(t_{1})^{-1}\ldots e(t_{k})^{-1}\sum_{0\leq
j_{1},\ldots,j_{k}\leq1}E_{j_{1}}(t_{1})\ldots
E_{j_{k}}(t_{k})D_{n}( x_{i_{1}}^{j_{1}}
\ldots x_{i_{k}}^{j_{k}})$$
where $e(\tau)=(1-s\tau)^{2},\,E_{0}(\tau)=1-s\tau ,\,E_{1}(\tau)=\tau.$
\end{theorem}
\noindent Some computations of $D$ polynomials are given at  the end of section 4. 

The expansion formula of Theorems \ref{expansion:formula} and \ref{expansion:D} reduce
the computations of link polynomials to
polynomials of positive braids with total exponent
$J=j_{1}+j_{2}+\ldots +j_{k}\leq k.$ If these braids (or some
positive conjugates of them) contain exponents greater than 1, apply
again and again expansion formula; finally we obtain \emph{simple
braids}. In  section 5 we give different characterizations  of this set of braids, the proofs and more properties are contained in \cite{Barbu-Rehana:02}. The first phrase in the next theorem could be taken as a definition :
\begin{theorem}\label{markov: th}a) A positive braid in $\mathcal{B}_{n}$ is simple if and only if it is conjugate to $$\beta_{A}=(x_{1}x_{2}\ldots x_{s_{1}-1})(x_{s_{1}+1}\ldots
x_{s_{2}-1})\ldots (x_{s_{r-1}+1}\ldots
x_{s_{r}-1})$$ for an decreasing  sequence $A=(a_{1},a_{2},\ldots , a_{r})$ with $a_{1}\geq a_{2}\geq \ldots \geq a_{r}\geq 2$ and $s_{i}=a_{1}+a_{2}+\ldots+a_{i}$.\\
b) HOMFLY polynomial of a closed simple braid $B_{A}$ is given by $$P_{n}(\beta_{A})(l,m)=\Big(-\frac{l+l^{-1}}{m}\Big)^{n-s_{r}+r-1}.$$
\end{theorem}
\noindent The values of Alexander-Conway, Jones and $D$ polynomials for simple braids are computed in section 5.

In the last section we show that Alexander-Conway, Jones and $D$ polynomials of $2$-braids are linearly independent over $\mathbb{C}[s,s^{-1}]$ and also
 \begin{theorem}\label{algebraically independent} The link polynomials   $\nabla$, $V$ and $D$ are algebraically independent  over $\mathbb{C}[s,s^{-1}].$
  \end{theorem}
\noindent The last example is a $2$-link with the same Alexander-Conway and Jones polynomials as the unlink, but with a nontrivial $D$ invariant.
  \medskip

\section{\textbf{ Nice Specializations}}

\medskip

\emph{Proof of Theorem \ref{recurrence: theorem}:} Consider the braid
(the exponent $e$ is positive)
$$\beta(e+2)=x_{i_{1}}^{a_{1}}
x_{i_{2}}^{a_{2}}\ldots x_{i_{j-1}}^{a_{j-1}}
x_{i_{j}}^{e+2}x_{i_{j+1}}^{a_{j+1}}\ldots x_{i_{k}}^{a_{k}}.$$

 Applying the skein relation for HOMFLY
 polynomial of this braid, the geometrical change appears
as a change at the last crossing of the j-th position. Both parts
$x_{i_{1}}^{a_{1}}x_{i_{2}}^{a_{2}}\ldots x_{i_{j-1}}^{a_{j-1}}$ and
$x_{i_{j+1}}^{a_{j+1}}\ldots x_{i_{k-1}}^{a_{k-1}}x_{i_{k}}^{a_{k}}$
remain unchanged. The local picture is:
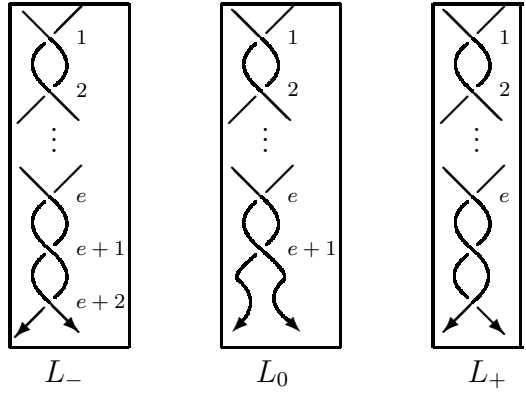
\begin{center}
\begin{picture}(360,160)
\thicklines
 \qbezier(80,140)(92,130)(82,122)      \qbezier(78,138)(68,130)(80,120)
 \put(80,140){\line(-1,1){10}}         \put(82,142){\line(1,1){11}}
 \put(80,120){\line(1,-1){11}}         \put(78,118){\line(-1,-1){10}}
 \put(90,138){{\tiny 1}}               \put(90,118){{\tiny 2}}
 \put(65,23){\line(0,1){130}}          \put(110,23){\line(0,1){130}}
 \put(145,23){\line(0,1){130}}         \put(190,23){\line(0,1){130}}
 \put(225,23){\line(0,1){130}}         \put(258,23){\line(0,1){130}}
 \put(65,23){\line(1,0){45}}           \put(145,23){\line(1,0){45}}
 \put(225,23){\line(1,0){34}}          \put(65,153){\line(1,0){45}}
 \put(145,153){\line(1,0){45}}         \put(225,153){\line(1,0){34}}
 \qbezier(160,140)(172,130)(162,122)   \qbezier(158,138)(148,130)(160,120)
 \put(160,140){\line(-1,1){11}}        \put(160,120){\line(1,-1){11}}
 \put(162,142){\line(1,1){10}}         \put(158,118){\line(-1,-1){10}}
 \put(170,138){{\tiny 1}}              \put(170,118){{\tiny 2}}
 \qbezier(240,140)(252,130)(242,122)   \qbezier(238,138)(228,130)(240,120)
 \put(240,140){\line(-1,1){11}}        \put(240,120){\line(1,-1){11}}
 \put(242,142){\line(1,1){10}}         \put(238,118){\line(-1,-1){10}}
 \put(250,138){{\tiny 1}}              \put(250,118){{\tiny 2}}
 \qbezier(80,80)(92,70)(82,62)         \qbezier(78,78)(68,70)(80,60)
 \qbezier(80,60)(92,50)(82,42)         \qbezier(78,58)(68,50)(80,40)
 \put(80,80){\line(-1,1){11}}          \put(82,82){\line(1,1){10}}
 \put(80,40){\vector(1,-1){11}}        \put(78,38){\vector(-1,-1){11}}
 \put(90,78){{\tiny $e$}}              \put(90,58){{\tiny $e+1$}}
 \put(90,38){{\tiny $e+2$}}            \put(160,80){\line(-1,1){11}}
 \put(162,82){\line(1,1){10}}          \qbezier(160,80)(172,70)(162,62)
 \qbezier(158,78)(148,70)(160,60)      \qbezier(160,60)(172,50)(168,48)
 \qbezier(168,48)(160,40)(172,32)      \put(172,32){\vector(1,-1){3}}
 \qbezier(158,58)(148,50)(152,48)      \qbezier(152,48)(160,40)(152,32)
 \put(152,32){\vector(-1,-1){3}}       \put(170,78){{\tiny $e$}}
 \put(170,58){{\tiny$e+1$}}            \put(240,80){\line(-1,1){11}} 
 \put(242,82){\line(1,1){10}}          \qbezier(240,80)(252,70)(242,62)
 \qbezier(238,78)(228,70)(240,60)      \qbezier(240,60)(252,50)(240,40)
 \qbezier(238,58)(228,50)(238,42)      \put(240,40){\vector(-1,-1){11}}
 \put(242,38){\vector(1,-1){10}}       \put(250,78){{\tiny $e$}}
 \put(78,10){$L_{-}$}                  \put(158,10){$L_{0}$}
 \put(238,10){$L_{+}$}                 \put(80,97){\vdots}  
 \put(160,97){\vdots}                  \put(240,97){\vdots}
\end{picture}
\end{center}
The skein relation $lP(L_{+})+l^{-1}P(L_{-})+mP(L_{0})=0$ with
correspondence $L_{-}\rightarrow \widehat{\beta}(e+2)$,
$L_{0}\rightarrow \widehat{\beta}(e+1),$ $L_{+}\rightarrow
\widehat{\beta}(e)$, gives
$$P_{n}(e+2)+mlP_{n}(e+1)+l^{2}P_{n}(e)=0.$$
The case $e<0$ can be reduced to the above case by  adding a new
factor on the $j+1$ position, $x_{i_{1}}^{a_{1}} \ldots
x_{i_{j}}^{e}\ldots x_{i_{k}}^{a_{k}}=x_{i_{1}}^{a_{1}} \ldots
x_{i_{j}}^{d+e}x_{i_{j}}^{-d}\ldots x_{i_{k}}^{a_{k}}$, with $d$
 big enough.
\begin{flushright}
$\Box$
\end{flushright}
\begin{example}\label{exp:homfly}A more complicated  recurrence is used  for the sequence $P(k)=P_{3}(\gamma_{k})$ where $\gamma_{k}$ is the 3-braid $x_{1}x_{2}x_{1}\ldots(k-\hbox{factors})$. This sequence  contains the powers of Garside braid $\Delta_{3}$: $\gamma_{6k}=\Delta_{3}^{2k},\,\gamma_{6k+3}=\Delta_{3}^{2k+1}.$

\begin{tabbing}
 0\=1\indent \= 2\qquad\qquad \=3 \=4\kill
\> 1) \>$P(2k+1)$\>=\>$ -mlP(2k)-l^{2}P(2k-1)$\\
\> 2)\> $P(6k+4)$\>=\>$-mlP(6k+3)-l^{2}P(6k+2)$\\
\> 3)\>$P(6k+2)$\>=\>$-mlP(6k+1)+ml^{3}P(6k-1)+l^{4}P(6k-2)$\\
\> 4)\> $P(6k)$\>=\>$-mlP(6k-1)+ml^{3}P(6k-3)-ml^{5}P(6k-5)-l^{6}P(6k-6).$
 \end{tabbing}
 The proof is the same as in \cite{Barbu-Nizami:09}. The characteristic polynomial of this recurrence has $0$ a triple root. In the case of Jones specialization the other three roots are $s^{6}$, $s^{12}$ and $s^{18}$. For Alexander-Conway polynomial, see section 3 (and for $D$ polynomial see section 4).
\end{example}
\quad The roots of the characteristic equation of the basic recurrence of Theorem \ref{recurrence: theorem}  are not
rational functions in $l,\,m$. From the point of view of recurrence
relations a "\emph{nice specialization}" $l\rightarrow
F(s),\,m\rightarrow G(s)$ will give "\emph{nice roots}" of the
characteristic equation.

\medskip

\begin{definition}
\label{rat:ional} Let us call a \emph{rational specialization} of
HOMFLY $P(l,m)$ a specialization $l=F(s),\,m=G(s)\,\, (F,G,$ are
rational functions$)$
such that:\\
\indent \emph{a)} $P  (\bigcirc\bigcirc)(F(s),G(s))$ is a Laurent polynomial in $s;$\\
\indent \emph{b)} the roots of the characteristic equation are $r_{1}= \lambda  s^{n},\,r_2= \mu s^{k}$
 where $n,\,k$ are   coprime integers and $\lambda ,\mu$ are non zero complex numbers.
\end{definition}
\begin{definition}
\label{nor:malized} A rational specialization is  \emph{
normalized} if $n$ equals 1 and \emph{nondegenerate} if  $n\neq k$.
\end{definition}
The next result says that there are few "nice specializations" and also that Laurent condition \emph{a)}  on the unlink gives Laurent polynomials for any link.
\begin{theorem}
\label{jones:conway} There are (up to  constants) only three
normalized  rational specialization of the HOMFLY polynomial, two nondegenerate and one degenerate:\\
\indent 1) Alexander-Conway with  roots $r_{1}=-s,\,r_{2}=s^{-1}$;\\
\indent 2) Jones with roots $r_{1}=-s,\,r_{2}=s^{3};$\\
\indent 3) roots $r_{1}=r_{2}=s.$
\end{theorem}
\begin{proof}The characteristic equation of the HOMFLY recurrence is
$$r^{2}+mlr+l^{2}=0.$$
Let $\lambda^{2}s^{n}$, $\mu^{2}s^{k}$ be the roots of the equation,
where
$\lambda^{2},\mu^{2}\in\mathbb{C^{\ast}},\,n,k\in\mathbb{Z},\,(n;k)=1.$\\
We should have $n+k=2q$ for some integer $q$, therefore
\begin{eqnarray*}
\label{eq:1} (\ast)\indent\indent\indent\indent l&=& \lambda\mu
s^{q},\,m=-\frac{\lambda}{\mu}s^{n-q}-\frac{\mu}{\lambda}s^{q-n}.
\end{eqnarray*}
Skein relation for eight figure $\infty$ gives:

\begin{eqnarray*}
  P(\bigcirc\bigcirc) &=& -\frac{l+l^{-1}}{m}=\frac{\lambda\mu s^{q}+(\lambda\mu s^{q})^{-1}}{\lambda\mu^{-1}s^{n-q}+\mu\lambda^{-1}s^{q-n}}\\
&=&\frac{\lambda^{2}\mu^{2}s^{2q}+1}{\lambda\mu s^q}\frac{\lambda\mu s^{n-q}}{\lambda^{2}s^{2(n-q)}+\mu^{2}}
=\frac{1}{s^{k}}\cdot\frac{\lambda^{2}\mu^{2}s^{2q}+1}{\lambda^{2}s^{2(n-q)}+\mu^{2}}
\end{eqnarray*}
and the first condition for a rational parametrization implies
$n-q\,|\,q$ or $n-q=0$. The case when $n-q=0$ gives the degenerate case $n=k=1,\,r_{1}=r_{2}=s$. If $n-q\,|\,q$, using $(n;k)=1$, we have $(n;q)=1$, hence $n-q=\pm1$. Normalization implies $n=1,$  hence $q=0,\,k=-1$ or
$q=2,\,k=3.$ One can find Alexander-Conway and Jones specializations (up to constants)
taking these values for $n$ and $q$ in $(\ast)$.
\end{proof}
Fibonacci recurrence for Jones polynomial is treated in \cite{Barbu-Nizami:09}. In the next two sections we  analyze the recurrence for Alexander-Conway  and for the degenerate case.

\section{\textbf{Alexander-Conway polynomial}}

\medskip
Using Alexander-Conway specialization, $l=\iota, m=\iota(s^{-1}-s),$
the recurrence relation for Alexander-Conway polynomials is
(Corollary \ref{cor:1})
$$\nabla_{n}(a_{1},\ldots,a_{j}+2,\ldots,a_{k})=(s^{-1}-s)\nabla_{n}(a_{1},\ldots,a_{j}+1,\ldots,a_{k})+\nabla_{n}(a_{1},\ldots,a_{j},\ldots,a_{k}).$$
%\medskip
The characteristic equation  $r^{2}=(s^{-1}-s)r+1$ has the roots
$r_{1}=-s,\,r_{2}=s^{-1}.$ Applying  computations with
multiple Fibonacci sequence (see \cite{nizami:08}, Theorem 3.8, for
complete details) one can find the expansion formula and compute
the generating function.

\emph{Proof of Theorem \ref{expansion:formula}:} Laurent polynomials
of (\cite{nizami:08}, Theorem 3.8 a) are given by:
\begin{eqnarray*}
  \Delta &=& r_{2}-r_{1}=\frac{s^{2}+1}{s} \\
  S_{0}^{[a]} &=& r_{1}^{a}r_{2}-r_{1}r_{2}^{a}=(-1)^{a}s^{a-1}+s^{1-a}=C^{a}(s) \\
  S_{1}^{[a]} &=& r_{2}^{a}-r_{1}^{a}=s^{-a}+(-1)^{a+1}s^{a}=C^{a+1}(s),
\end{eqnarray*}
therefore the general term is given by
$$\nabla_{n}(x_{i_{1}}^{a_{1}} \ldots
x_{i_{k}}^{a_{k}})(s)=\Big(\frac{s}{s^2+1}\Big)^{k}\sum_{0\leq
j_{1},\ldots,j_{k}\leq1}C^{a_{1}+j_{1}}(s)\ldots
C^{a_{k}+j_{k}}(s)\nabla_{n}( x_{i_{1}}^{j_{1}} \ldots
x_{i_{k}}^{j_{k}})(s).$$
\begin{flushright}
$\Box$
\end{flushright}

\emph{Proof of Theorem \ref{generating:function}}\indent From the second
part of (\cite{nizami:08}, Theorem 3.8), we use the polynomials of
formal variables of $t_{1},t_{2},\ldots,t_{k}$ (with coefficients in
the ring of Laurent polynomials in variable $s$):
\begin{eqnarray*}
  q(\tau) &=& (1-r_{1}\tau)(1-r_{2}\tau)=1+(s-s^{-1})\tau-\tau^{2}, \\
  Q_{0}(\tau) &=& 1-(r_{1}+r_{2})\tau=1+(s-s^{-1})\tau, \\
  Q_{1}(\tau) &=& \tau
\end{eqnarray*}
\noindent  and we obtain
$$\mathcal{AC}_{n}(t_{1},\ldots,t_{k})=q(t_{1})^{-1}\ldots q(t_{k})^{-1}\sum_{0\leq
j_{1},\ldots,j_{k}\leq1}Q_{j_{1}}(t_{1})\ldots
Q_{j_{k}}(t_{k})\nabla_{n}( x_{i_{1}}^{j_{1}}
\ldots x_{i_{k}}^{j_{k}}).$$

\begin{flushright}
$\Box$
\end{flushright}
Using the recurrence relation  we give some examples of Alexander-Conway  polynomials:
\begin{example}\label{C exp}The recurrence relation for $\nabla_{2}(x_{1}^{a})$ starts with $\nabla_{2}(1)=0,\,\nabla_{2}(x_{1})=1$ and gives the general term,
       $$ \nabla_{2}(x_{1}^{a}) =\frac{s^{-a}+(-1)^{a+1}s^{a}}{s+s^{-1}} = s^{1-a}-s^{3-a}+s^{5-a}-\ldots +(-1)^{a+1}s^{a-1},\,a\geq 0,$$
       and   $\nabla_{2}(x_{1}^{a}) =(-1)^{a+1}s^{a+1}+(-1)^{a}s^{a+3}+\ldots-s^{-a-3}+s^{-a-1},\,a<0$.
   \end{example}
   \begin{example}\label{AC:exp}With  specializations $(l,m)\longrightarrow (\iota,\iota(s^{-1}-s))$ in Example \ref{exp:homfly}, we have recurrences for $\nabla(k)=\nabla_{3}(\gamma_{k})$, $\gamma_{k}$ is the $3$-braid $x_{1}x_{2}x_{1}\ldots (k\hbox{-factors})$:
\begin{tabbing}
 0\=1\indent \= 2\qquad\qquad \=3 \=4\kill
\> 1)\> $\nabla(2k+1)$\>=\>$(s^{-1}-s)\nabla(2k)+\nabla(2k-1)$   \\
\> 2)\>  $\nabla(6k)$\> = \>$(s^{-1}-s)[\nabla(6k-1)+\nabla(6k-3)+\nabla(6k-5)]+\nabla(6k-6)$\\
\> 3)\>  $\nabla(6k+2)$\>=\>$(s^{-1}-s)[\nabla(6k+1)+\nabla(6k-1)]+\nabla(6k-2)$ \\
\> 4)\>  $\nabla(6k+4)$\>=\>$(s^{-1}-s)\nabla(6k+3)+\nabla(6k+2).$
\end{tabbing}

 \end{example}
\begin{remark} From these  recurrences we can compute Jordan normal form of  the recurrence matrix $M(k)$:
\medskip

$M(k)\sim\left(
           \begin{array}{llll}
             J_{3} & 0 & 0 & 0 \\
             0 & 1 & 0 & 0 \\
             0 & 0 & s^{6} & 0 \\
             0 & 0 & 0 & s^{-6} \\
           \end{array}
         \right)
$, where $J_{3}=\left(
            \begin{array}{ccc}
              0 & 1 & 0 \\
              0 & 0 & 1 \\
              0 & 0 & 0 \\
            \end{array}
          \right).
$
\end{remark}
\begin{example}
\label{first:9}The first six Alexander-Conway polynomials are\\
 \indent$\nabla(0)=0$ \\
  \indent $\nabla(1)=0 $\\
   \indent$\nabla(2) =1$\\
  \indent $\nabla(3)=s^{-1}-s$ \\
  \indent $\nabla(4)=s^{-2}-1+s^{2}$ \\
   \indent$\nabla(5)=s^{-3}-s^{-1}+s-s^{3}$ \\
\end{example}Using recurrences for  $\nabla(k)$ in Example \ref{AC:exp}, we have explicit formulae:
\begin{proposition}\label{B3:exp}
  \indent $\nabla(6k)=(s^{-4}-s^{-2})\sum^{k-1}\limits_{i=0}s^{-6i}+(s^{4}-s^{2})\sum^{k-1}\limits_{i=0}s^{6i} $\\
  \indent $\nabla(6k+1)=(s^{-5}-s^{-3})\sum^{k-1}\limits_{i=1}s^{-6i}-(s^{5}-s^{3})\sum^{k-1}\limits_{i=0}s^{6i}$\\
  \indent $\nabla(6k+2)=(s^{-6}-s^{-4})\sum^{k-1}\limits_{i=0}s^{-6i}+1+(s^{6}-s^{4})\sum^{k-1}\limits_{i=0}s^{6i}$\\
  \indent $\nabla(6k+3)=(s^{-7}-s^{-5})\sum^{k-1}\limits_{i=0}s^{-6i}+s^{-1}-s-(s^{7}-s^{5})\sum^{k-1}\limits_{i=0}s^{6i}$\\
  \indent $\nabla(6k+4)=(s^{-2}-1)\sum^{k}\limits_{i=0}s^{-6i}+1+(s^{2}-1)\sum^{k}\limits_{i=0}s^{6i}$\\
  \indent $\nabla(6k+5)=(s^{-3}-s^{-1})\sum^{k}\limits_{i=0}s^{-6i}-(s^{3}-s^{1})\sum^{k}\limits_{i=0}s^{6i}$\\
In the first four relations $k$ is positive, in the last two, $k$ is nonnegative.
\end{proposition}
\begin{definition}For a non zero Laurent polynomial
 $P=a_{p}s^{p}+a_{p-1}s^{p-1}+\ldots +a_{q+1}s^{q+1}+a_{q}s^{q} \in \mathbb{C}[s,s^{-1}]$, where $p$ and $q$ are integers and $a_{p}\neq0,\,a_{q}\neq0,\,p\geq q$, $p$ is called the \emph{degree of} $P,\,\deg(P)$, $q$ is
called the \emph{order  of} $P,\,\hbox{ord}(P)$ and $p-q+1$ is called the \emph{breadth of} $P,\,\br(P)$. For $p=0$ we take $\hbox{deg}(0)=-\infty, \br(0)=0$.
\end{definition}

\begin{proposition}
\label{deg:inequ}
If $a_i \in  \{0,1\},$ then
\end{proposition}
\vspace{-0.4cm}
%\nabla_{3}(\underbrace{x_1x_2x_1x_2\ldots x_1x_2}_{2k \,factors})
 $$\hbox{\deg}\nabla_{3}(\underbrace{x_1^{a_1}x_2^{a_2}\ldots
x_1^{a_{2k-1}}x_2^{a_{2k}}}_{2k\, factors})\leq 2k-2.$$

 \begin{proof}We prove the result by induction: for $k=1$,
$$\nabla(x_{1}^{0}x_{2}^{0})=\nabla(x_{1}^{1}x_{2}^{0})=\nabla(x_{1}^{0}x_{2}^{1})=0,\,\nabla(x_{1}x_{2})=1.$$
Suppose that the result is true for $k-1$, i.e $\deg\nabla_{3}(x_{1}^{a_{1}}x_2^{a_{2}}\ldots
x_{1}^{a_{2k-3}}x_2^{a_{2k-2}})\leq 2k-4.$\\

For $\beta=x_1^{a_{1}}x_2^{a_{2}}\ldots x_{2}^{a_{2k-2}}x_{1}^{a_{2k-1}}x_{2}^{a_{2k}}$, we have two cases:\\

\emph{Case I:} If $a_{1}=a_{2}=\ldots =a_{2k}=1$, we have equality (see Proposition \ref{B3:exp}).\\

\emph{Case II:} If there is at least one zero, say $a_{2k}=0$, then
$\beta$ is conjugate to $x_{1}^{a_{1}+a_{2k-1}}x_{2}^{a_{2}}\ldots x_{1}^{a_{2k-3}}x_{2}^{a_{2k-2}}$.
If $a'_{1}=a_{1}+a_{2k-1}\leq1$, then it is clear using induction step, otherwise use expansion formula\\
\begin{eqnarray*}
\nabla_{3}(x_1^{a_1}x_2^{a_2}\,.\,.\,
x_1^{a_{2k-1}}x_2^{a_{2k}})&=& \nabla_{3}(x_1^{2}x_2^{a_2}\ldots
x_2^{a_{2k-2}})\\
&=& \Big(\frac{s}{s^2+1}\Big )^{2k-2}\sum _{0\leq
j_{*}\leq1}C^{2+j_{1}}\,.\,.\,
C^{a_{2k-2}+j_{2k-2}}\nabla_{3}( x_{1}^{j_{1}}\,.\,.\,
x_{2}^{j_{2k-2}})
\end{eqnarray*}

$\deg(\frac{s}{s^2+1})^{2k-2}=-(2k-2)$. As $a_{i}+j_{i}\in \{0,1,2\}$ for $i=2,\ldots,2k-2$ and $C^{0}=C^{2}=s+s^{-1}$, $C^{1}=0$, $C^{3}=-s^{2}+s^{-2}$

 $ \hbox{max}\{\deg(C^{2+j_{1}}\ldots
C^{a_{2k-2}+j_{2k-2}}),\,\hbox{where} \,\,0\leq j_{1},j_{2},\ldots,j_{2k-2}\leq 1\}$\\
$\indent =2+\hbox{max}\{\deg( C^{a_{2}+j_{2}}\ldots
C^{a_{2k-2}+j_{2k-2}}),\,\hbox{where} \,\,0\leq j_{2},\ldots,j_{2k-2}\leq 1\}$\\
$\indent =2+2k-3=2k-1$\\
Also from induction hypothesis $\deg \nabla_{3}( x_{1}^{j_{1}}\ldots x_{2}^{j_{2k-2}})\leq 2k-4$. So
$$\deg\nabla_{3}(x_1^{a_1}x_2^{a_2}\ldots
x_1^{a_{2k-1}}x_2^{a_{2k}})\leq-(2k-2)+2k-1+2k-4\leq2k-2.$$
\end{proof}
For large exponents, we have
\begin{theorem}
If all $a_i\geq2,$ then
\end{theorem}
\vspace{-0.4cm}
 $$\hbox{\deg} \nabla_{3}(x_1^{a_1}x_2^{a_2}\ldots
x_1^{a_{2k-1}}x_2^{a_{2k}})=a_1+a_2+\ldots+a_{2k}-2.$$

\begin{proof}To evaluate $\nabla_{3}(x_1^{a_1}x_2^{a_2}\ldots
x_1^{a_{2k-1}}x_2^{a_{2k}})$, use the expansion formula and we
have
\begin{eqnarray*}
\nabla_{3}(x_1^{a_1}x_2^{a_2}\,.\,.\,
x_1^{a_{2k-1}}x_2^{a_{2k}})&=& \Big(\frac{s}{s^2+1}\Big )^{2k}\sum _{0\leq
j_{*}\leq1}C^{a_{1}+j_{1}}\,.\,.\,
C^{a_{2k}+j_{2k}}\nabla_{3}( x_{1}^{j_{1}}\,.\,.\,
x_{2}^{j_{2k}})
\end{eqnarray*}
$\deg(\frac{s}{s^{2}+1})^{2k}=-2k,$ Proposition \ref{deg:inequ} and Proposition \ref{B3:exp} imply that
$$\deg \nabla_{3}(\underbrace{x_1x_2x_1x_2\ldots x_1x_2}_{2k
\,\hbox{factors}})=2k-2\,\,\,\,\, (j_{1}=j_{2}=\ldots j_{2k}=1)$$\\
 As  $\deg C^{a_{i}+1}=a_{i}>a_{i}-1=\deg C^{a_{i}+0}$, so\\
$\max\{\deg(C^{a_{1}+j_{1}}C^{a_{2}+j_{2}}\ldots C^{a_{2k}+j_{2k}}),\, \hbox{where} \,\,0\leq j_{1},j_{2},\ldots,j_{2k}\leq 1\}$\\
$=\deg(C^{a_{1}+1}C^{a_{2}+1}\ldots C^{a_{2k}+1})=a_{1}+a_{2}+\ldots +a_{2k}.$\\
Finally\begin{eqnarray*}
 \deg \nabla_{3}(x_1^{a_1}x_2^{a_2}\ldots
x_1^{a_{2k-1}}x_2^{a_{2k}})&=&-2k+a_{1}+a_{2}+\ldots+a_{2k}+2k-2\\
&=&a_{1}+a_{2}+\ldots +a_{2k}-2.
\end{eqnarray*}

\end{proof}

\begin{cor}If all $a_i\geq2,$ then
$$\br( \nabla_{3}(x_1^{a_1}x_2^{a_2}\ldots
x_1^{a_{2k-1}}x_2^{a_{2k}}))=2(a_1+a_2+\ldots+a_{2k})-3.$$
\end{cor}

\medskip

\section{\textbf{The degenerate case}}
The specialization $l=s,m=-2$, corresponding to the roots $r_{1}=r_{2}=s$, gives \emph{the degenerate polynomial $D$} with skein relation $$sD(L_{+})+s^{-1}D(L_{-})-2D(L_{0})=0$$ and recurrence relation $$D_{n}(a_{1},\ldots, a_{j}+2, \ldots, a_{k})=2sD_{n}(a_{1},\ldots, a_{j}+1, \ldots, a_{k})-s^{2}D(a_{1},\ldots, a_{j}, \ldots, a_{k}).$$

  \medskip

\emph{Proof of Theorem \ref{expansion:D}}: As in \cite{Barbu-Nizami:09} for Jones polynomial and as in section 3 for Alexander-Conway polynomial we find the basic polynomials
  \begin{eqnarray*}
     R_{0}^{[a]}(s)&=& -r_{1}^{a-1}r_{2}-r_{1}^{a-2}r_{2}^{2}-\ldots -r_{1}r_{2}^{a-1}=(1-a)s^{a}\\
    R_{1}^{[a]}(s) &=& r_{1}^{a-1}+r_{1}^{a-2}r_{2}-\ldots +r_{2}^{a-1}=as^{a-1}
  \end{eqnarray*}
  Using the notations and computations of  \cite{nizami:08}, Remark 2.2,  we obtain  the result.

\begin{flushright}
$\Box$
\end{flushright}

For a fixed sequence $(x_{i_{1}}x_{i_{2}}\ldots x_{i_{k}})$ in $\mathcal{B}_{n}$ the \emph{generating function} for the degenerate polynomial is defined by the formal series in $t_{1},\ldots ,t_{k}:$
$$\mathcal{D}_{n}(t_{1},\ldots ,t_{k})=\sum_{(a_{1},a_{2},\ldots,a_{k})\in\mathbb{Z}^{k}}D_{n}(x_{i_{1}}^{a_{1}} \ldots x_{i_{k}}^{a_{k}})(s)t_{1}^{a_{1}}\ldots t_{k}^{a_{k}}.$$
\emph{Proof of Theorem \ref{generating:D}:} Using \cite{nizami:08}, Theorem 3.8,
the corresponding polynomials are $e(\tau)=(1-s\tau)^{2}$ and $E_{0}(\tau)=1-2s\tau ,\,E_{1}(\tau)=\tau.$

\begin{flushright}
$\Box$
\end{flushright}

\begin{proposition}(\textbf{Properties of $D$ polynomial})\label{properties:D} Let $L$, $L_{1}$, $L_{2}$ be oriented links,  $\overline{L}$ be the mirror image of $L$, $\overleftarrow{L}$ the link with reversed orientation of every component of $\overrightarrow{L}$, $"\small\coprod"$ be the distant union of two links, then

\medskip

\indent \emph{1)} $D(L_{1}\coprod L_{2})=\frac{1}{2s}(s^{2}+1)D(L_{1})D(L_{2})$\\
\indent \emph{2)} $D(\overline{L})(s)=D(L)(s^{-1})$\\
\indent \emph{3)}  $D({L_{1}}+L_{2})=D(L_{1})D(L_{2})$\\
\indent \emph{4)} $D(\overleftarrow{L})=D(\overrightarrow{L})$.
\end{proposition}
\begin{proof}The proofs are as in \cite{Lickorish:97} for Jones polynomial.
\end{proof}
\begin{remark}a) For any link $L$, $D(L)(1)=1$.\\
b) The degenerate polynomials take values in $\mathbb{Z}[\frac{1}{2}][s,s^{-1}]$.
\end{remark}
\begin{example}\label{D exp}Starting with $D(\bigcirc)=D_{2}(x_{1})=1$ we obtain $D(\bigcirc\bigcirc)=D_{2}(1)=\frac{1}{2}(s+s^{-1})$ and $D_{2}(x_{1}^{n})=\frac{1}{2}[(1-n)s^{n+1}+(1+n)s^{n-1}]$; in particular  $D_{2}(x_{1}^{2})=\frac{1}{2}(-s^{3}+3s)$ and $D_{2}(x_{1}^{3})=-s^{4}+2s^{2}$.
\end{example}
\begin{example}With  specializations  $(l,m)\longrightarrow (s,-2)$ in Example \ref{exp:homfly}, we have the following recurrences for $D(k)=D_{3}(\gamma_{k})$ ($\gamma_{k}$ is the 3-braid $x_{1}x_{2}x_{1}\ldots \,\hbox{with}\,\,k-\hbox{factors}$):

\begin{tabbing}
  0\=1\indent \= 2\qquad\qquad \=3 \=4\kill
  % \> for next tab, \\ for new line...
  \>1) \> $D(2k+1)$\> =\>  $2sD(2k)-s^{2}D(2k-1)$ \\
  \>2) \> $D(6k+4)$\>=\>$2sD(6k+3)-s^{2}D(6k+2)$\\
  \>3) \> $D(6k+2)$ \>=\> $2sD(6k+1)-2s^{3}D(6k-1)+s^{4}D(6k-2)$\\
  \>4) \> $ D(6k)$ \>=\> $2sD(6k-1)-2s^{3}D(6k-3)+2s^{5}D(6k-5)-s^{6}D(6k-6).$\\
 \end{tabbing}
  \end{example}
  \begin{remark} These recurrences relations for $D$ polynomials give a recurrence matrix with Jordan normal form:

\medskip
$M(k)\sim \left(
           \begin{array}{llll}
             J_{3} & 0 & 0 & 0 \\
             0 & s^{6} & 0 & 0 \\
             0 & 0 & \gamma s^{6} & 0 \\
             0 & 0 & 0 & \overline{\gamma} s^{6} \\
           \end{array}
         \right)
$, with $J_{3}=\left(
            \begin{array}{ccc}
              0 & 1 & 0 \\
              0 & 0 & 1 \\
              0 & 0 & 0 \\
            \end{array}
          \right)
$, $\gamma=4+\iota  \sqrt{3}$.
\end{remark}
\begin{example}\label{DD exp}First few $D(k)=D_{3}(\gamma_{k})$ polynomials  are:\\
$D(0)$=$\frac{1}{4}(s^{2}+2+s^{-2})$, $D(1)=\frac{1}{2}(s+s^{-1})$, $D(2)=1$,  $D(3)=\frac{1}{2}(-s^{3}+3s)$,\\
\medskip
\noindent $D(4)=-s^{4}+2s^{2}$, $D(5)=\frac{1}{2}(-3s^{5}+5s^{3})$, $D(6)=\frac{1}{4}(-s^{8}-6s^{6}+11)$ and for $n\geq7,D(n)=\frac{1}{4}(-a_{n}s^{n+2}+b_{n}s^{n+1}+c_{n}s^{n})$,  where the sequence $(a_{n})_{n\geq6}$ is given by $a_{12k}=12k-1$, $a_{12k+i}=12k$ for $i=1,\ldots,5$, $a_{12k+6}=12k+1$, $a_{12k+6+i}=12k+2i$ for $i=1,\ldots,5$.
\end{example}

\begin{example}If $\gamma=\alpha x_{n}^{k}\beta$ is an $n+1$-braid with $\alpha,\beta\in\mathcal{B}_{n}$, then
\medskip

$D_{n+1}(\gamma)=D_{n}(\alpha\beta)D_{2}(x_{1}^{k}).$
In particular $D_{3}(x_{1}^{a}x_{2}^{b})=D_{2}(x_{1}^{a})D_{2}(x_{1}^{b}).$
\end{example}

 \begin{proposition}\label{D deg}If $a_{i}\geq 0$ and the $3$-braid $x_{1}^{a_{1}}\ldots x_{2}^{a_{2k}}$ has degree $A=\sum^{2k}\limits_{i=1}a_{i}$ then $$\deg D_{3}(x_{1}^{a_{1}}x_{2}^{a_{2}}\ldots x_{1}^{a_{2k-1}}x_{2}^{a_{2k}})\leq A+2.$$
 \end{proposition}
 \begin{proof}By induction on $k$ and by previous example, this inequality is true for exponents $0\leq a_{i}\leq1$. For arbitrary non negative exponents the general expansion formula gives\\
 $$D_{3}(a_{1},\ldots,a_{2k})=\sum_{J_{*}\in \{0,1\}^{2k}}R_{J_{*}}^{A_{*}}D_{3}(j_{1},j_{2},\ldots,j_{2k}),$$ where $R_{J_{*}}^{A_{*}}=R_{j_{1}}^{[a_{1}]}\ldots R_{j_{k}}^{[a_{2k}]}$ has the degree at most $\sum^{2k}\limits_{i=1}(a_{i}-j_{i})=A-J$ (zero occurs only if $(a_{i},j_{i})=(0,1)\,\hbox{or}\,(1,0)$) and $D_{3}(j_{1},\ldots,j_{2k})$ has degree $\leq J+2.$
 \end{proof}

 \begin{example}$D_{3}(x_{1}^{a_{1}}x_{2}^{a_{2}})=\frac{1}{4}(1-a_{1})(1-a_{2})s^{A+2}+\ldots$
  \end{example}

 \begin{example}$D_{3}(x_{1}^{a_{1}}x_{2}^{a_{2}}x_{1}^{a_{3}}x_{2}^{a_{4}})=\frac{1}{4}[1-A+(a_{1}+a_{3})(a_{2}+a_{4})-a_{1}a_{2}a_{3}a_{4}]s^{A+2}+\ldots$
  \end{example}
\begin{remark}In these cases ($k=1$ and $k=2$), for exponents $a_{i}\geq2$, we have equality in Proposition \ref{D deg}. Perhaps the same is true for arbitrary $k$.
\end{remark}

\section{\textbf{HOMFLY polynomial for simple braids}}
\medskip
We describe a subset of divisors of Garside braid $\Delta_{n}$ which is related with HOMFLY recurrence and also has nice properties with respect to conjugation in $\mathcal{B}_{n}.$ We denote by $\mathcal{MB}_{n}$ the monoid of positive braids and by $\Div(\Delta_{n})$ the set of subwords of $\Delta_{n}$ (which coincides with the set of left divisors, the set of right divisors of $\Delta_{n}$, and also with the set of square free elements of  $\mathcal{MB}_{n}$, see \cite{Garside:69}). It is also well known that conjugation of positive braids $\alpha$, $\alpha'$ in $\mathcal{B}_{n}$ is equivalent with conjugation in $\mathcal{MB}_{n}$: there is a positive braid $\beta$ such that $\alpha\beta=\beta\alpha'$. Initial terms in the recurrence relation for HOMFLY are polynomials $P(\widehat{\beta})=P_{n}(\beta)$ for positive braids which have, up to a conjugation, all exponents equal to 1. One can define this subset, denoted by $\mathcal{SB}_{n}$, of positive braids in many ways:
\begin{definition}\label{simple braid} $\beta \in \mathcal{MB}_{n}$ is \emph{simple} if one presentation of $\beta$  contains no repeated generators $x_{i}$ (in this case, any presentation of  $\beta$ as a positive braid  has the same property).
\end{definition}

\begin{definition}$\beta\in \mathcal{MB}_{n}$ is \emph{simple} if all its conjugates in $\mathcal{MB}_{n}$ are square free.
\end{definition}
The equivalence of these definitions is proved in \cite{Barbu-Rehana:02}. In the same paper one can find canonical forms for simple braids as a product of disjoint cycles and the canonical form for conjugacy class of a simple braid (see Theorem \ref{markov: th} a)). Here is a picture of a simple braid: if $A=(4,3,2,2)$ and $\beta_{A}=(x_{1}x_{2}x_{3})(x_{5}x_{6})(x_{8})(x_{10})$ is  in $\mathcal{B}_{13}$, then
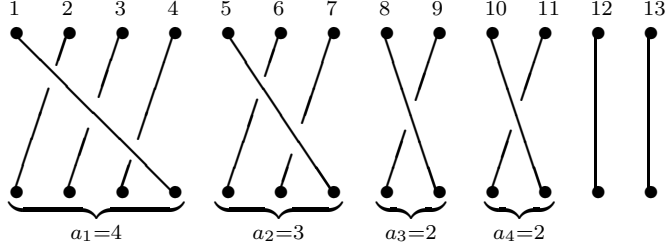
\begin{center}
\begin{picture}(360,90)
\thicklines
   \put(40,80){{\tiny 1}}      \put(40,10){$\bullet$}
   \put(40,70){$\bullet$}      \put(42,73){\line(1,-1){60}}
   \put(42,13){\line(1,3){13}} \put(60,80){{\tiny 2}}
   \put(60,10){$\bullet$}      \put(60,70){$\bullet$}
   \put(62,13){\line(1,3){8}}  \put(62,73){\line(-1,-3){4}}
   \put(80,80){{\tiny 3}}      \put(80,10){$\bullet$}
   \put(80,70){$\bullet$}      \put(82,73){\line(-1,-3){8}}
   \put(82,13){\line(1,3){4}}  \put(100,80){{\tiny 4}}
   \put(100,10){$\bullet$}     \put(100,70){$\bullet$}
   \put(102,72){\line(-1,-3){13}} \put(120,80){{\tiny 5}}
   \put(120,10){$\bullet$}     \put(120,70){$\bullet$}
   \put(123,72){\line(2,-3){40}} \put(122,12){\line(1,3){12}}
   \put(140,80){{\tiny 6}}     \put(140,10){$\bullet$}
   \put(140,70){$\bullet$}     \put(142,72){\line(-1,-3){5}}
   \put(142,12){\line(1,3){5}} \put(160,80){{\tiny 7}}
   \put(160,10){$\bullet$}     \put(160,70){$\bullet$}
   \put(162,72){\line(-1,-3){11}} \put(180,80){{\tiny 8}}
   \put(180,10){$\bullet$}     \put(180,70){$\bullet$}
   \put(200,80){{\tiny 9}}     \put(200,10){$\bullet$}
   \put(200,70){$\bullet$}     \put(182,72){\line(1,-3){20}}
   \put(182,12){\line(1,3){8}} \put(202,72){\line(-1,-3){8}}
   \put(220,80){{\tiny 10}}    \put(220,10){$\bullet$}
   \put(220,70){$\bullet$}     \put(240,80){{\tiny 11}}
   \put(240,10){$\bullet$}     \put(240,70){$\bullet$}
   \put(222,72){\line(1,-3){20}} \put(222,12){\line(1,3){8}}
   \put(242,72){\line(-1,-3){8}} \put(260,80){{\tiny 12}}
   \put(260,10){$\bullet$}     \put(260,70){$\bullet$}
   \put(262,72){\line(0,-1){60}} \put(280,80){{\tiny 13}}
   \put(280,10){$\bullet$}     \put(280,70){$\bullet$}
   \put(282,72){\line(0,-1){60}}  \put(40,10){{$\underbrace{\,\,\,\,\,\,\,\,\,\,\,\,\,\,\,\,\,\,\,\,\,\,\,\,\,\,\,\,\,\,\,\,\,}_{a_{1}=4}$}}
   \put(118,10){{$\underbrace{\,\,\,\,\,\,\,\,\,\,\,\,\,\,\,\,\,\,\,\,\,\,\,\,}_{a_{2}=3}$}}
   \put(179,10){{$\underbrace{\,\,\,\,\,\,\,\,\,\,\,\,\,}_{a_{3}=2}$}}
   \put(220,10){{$\underbrace{\,\,\,\,\,\,\,\,\,\,\,\,\,}_{a_{4}=2}$}}
\end{picture}
\end{center}

\bigskip

It is obvious that the closed braid $\widehat{\beta_{A}}$ is the trivial link with $n-s_{r}+r=n-\deg(\beta)$ components.
\begin{theorem}The value of  HOMFLY polynomial $P_{n}(\beta)$ of a simple braid $\beta$ is:
$$P_{n}(\beta)=\Big(-\frac{l+l^{-1}}{m}\Big)^{n-deg(\beta)-1}.$$
\end{theorem}
\begin{proof}The closure of the simple braid  $\beta$ is the trivial link with  $n-\deg(\beta)$ components and $P_{n}(\hbox{unlink with $k$ components})=\Big(-\frac{l+l^{-1}}{m}\Big)^{k-1}.$ Alternatively, one can use Markov moves $d=\deg(\beta)$ times and obtain $P_{n}(\beta)=P_{n-d}(1)$.
\end{proof}
\noindent\emph{Proof of Theorem \ref{markov: th} b):} The degree of the simple braid
 $$\beta_{A}=(x_{1}x_{2}\ldots x_{s_{1}-1})(x_{s_{1}+1}\ldots
x_{s_{2}-1})\ldots (x_{s_{r-1}+1}\ldots
x_{s_{r}-1})$$ is $s_{r}-r$.
\begin{flushright}
$\Box$
\end{flushright}

With  specializations $(l,m)\longrightarrow (\iota,\iota(s^{-1}-s)),(l,m)\longrightarrow (\iota s^{2},\iota(s-s^{-1}))$ and $(l,m)\longrightarrow (s,-2)$ respectively, we obtain the initial values for recurrence:

\begin{cor}$\nabla_{n}(\beta_{A})=\left\{
  \begin{array}{ll}
     & 0, \, \hbox{if $\deg(\beta)\neq n-1$}\\
     & 1, \, \hbox{if $\deg(\beta)= n-1$}.
  \end{array}
\right.$

\end{cor}
\begin{cor}\label{V cor}
$V_{n}(\beta_{A})=(-\frac{s^{2}+1}{s})^{n-deg(\beta)-1}.$
\end{cor}

\begin{cor}\label{Cor D}
$D_{n}(\beta_{A})=\Big(\frac{s^{2}+1}{2s}\Big)^{n-deg(\beta)-1}.$
\end{cor}

\section{\textbf{Relation between link polynomials}}
Using previous computations for some families of braids (Examples \ref{C exp}, \ref{D exp}, \ref{DD exp} and Corollaries \ref{V cor} and \ref{Cor D}) we can give "theoretical " proofs of the independence of the link polynomials $\nabla$, $V$ and $D$. More important for knot theorist is an example where $D$ can differentiate links and one of the links of \cite{trivial: Jones} is good for this.
\begin{definition} A set of $k$ families of Laurent polynomials  $\{f_{1}^{\alpha},f_{2}^{\alpha},\ldots,f_{k}^{\alpha}\}_{\alpha \in\mathcal{A}},\,\,f_{i}^{\alpha}\in\mathbb{C}[s,s^{-1}]$ is called \emph{linearly independent } over $\mathbb{C}[s,s^{-1}]$ if from an equation $$\sum^{k}\limits_{i=1}g_{i}(s,s^{-1})f_{i}^{\alpha}=0,\,\,\,(g_{i}(s,s^{-1})\in\mathbb{C}[s,s^{-1}])$$
which is satisfied for any value of parameter $\alpha$ in $\mathcal{A}$ we obtain $g_{i}=0,\,i=1,\ldots,k.$
 \end{definition}
 \begin{definition}A set of $k$ families of Laurent polynomials $\{f_{1}^{\alpha},\ldots,f_{k}^{\alpha}\}_{\alpha\in\mathcal{A}}$, $f_{i}^{\alpha}\in\mathbb{C}[s,s^{-1}]$ are called \emph{algebraically independent} over $\mathbb{C}[s,s^{-1}]$ if from a polynomial equation  $\Phi(f_{1}^{\alpha},\ldots,f_{k}^{\alpha})=0$, with $\Phi\in\mathbb{C}[s,s^{-1}][T_{1},\ldots,T_{k}]$,  which is satisfied for any $\alpha\in\mathcal{A}$, we obtain $\Phi=0$.
 \end{definition}
 As a first example we take $\mathcal{A}=\mathcal{B}_{2}$, $f_{1}(\beta)=\nabla_{2}(\beta)$, $f_{2}(\beta)=V_{2}(\beta)$, $f_{3}(\beta)=D_{2}(\beta)$ and we obtain
 \begin{proposition}$\{\nabla_{2}(\beta),V_{2}(\beta),D_{2}(\beta)\}_{\beta\in\mathcal{B}_{2}}$ is a linearly independent family over $\mathbb{C}[s,s^{-1}]$.
 \end{proposition}
 \begin{proof}Suppose we have  relations $A\nabla_{2}(x_{1}^{n})+BV_{2}(x_{1}^{n})+CD_{2}(x_{1}^{n})=0$ where $A,B,C\in\mathbb{C}[s,s^{-1}]$ and $\deg V_{2}(x_{1}^{n})=3n-1$, $\deg \nabla_{2}(x_{1}^{n})=n-1 $, $\deg D_{2}(x_{1}^{n})=n+1\,(n\geq2)$. If $B$ in not zero, take $b=\deg B$ and choose $n$ big enough to have $b+(3n-1)>\hbox{max}\{\deg A,\deg C\}+n+1$. Then the equation  $$A(s,s^{-1})V_{2}(x_{1}^{n})+B(s,s^{-1})\nabla_{2}(x_{1}^{n})+C(s,s^{-1})D_{2}(x_{1}^{n})=0$$ gives a contradiction. Therefore $B=0$ and for $A$, $C$ non zero we find a contradiction comparing the order of $\nabla_{2}(x_{1}^{n})$ (which is $1-n$)  and the order of $D_{2}(x_{1}^{n})$, which is positive.
 \end{proof}

 In the second example we choose $\mathcal{A}=\mathcal{SB}=\coprod\limits_{n}\mathcal{SB}_{n}$ the set of positive braids, $f_{1}(\beta)=V(\widehat{\beta})$ and $f_{2}(\beta)=D(\widehat{\beta})$.

\begin{proposition}$\{V(\widehat{\beta}),D(\widehat{\beta})\}_{\beta\in\mathcal{SB}}$ is an algebraic independent family over $\mathbb{C}[s,s^{-1}].$
\end{proposition}
\begin{proof}
Take $P(Y, Z)$ a polynomial with coefficients in $\mathbb{C}[s, s^{-1}]$ which is zero for $Y=V(\widehat{\beta}), Z=D(\widehat{\beta})$ and any simple braid $\beta$. We will show that $P$ is identically zero. Looking at the family  $\beta^{n}_{k}=x_{1}x_{2} \ldots x_{k-1}, \,\ 2\leq k\leq n-1$ and $k \equiv n \mod(2),$ we find $V(\widehat{\beta^{n}_{k}})=(\frac{s^2+1}{s})^{n-k},$ and $D(\widehat{\beta^{n}_{k}})=(\frac{s^2+1}{2s})^{n-k}$ and also $P(V(\widehat{\beta^{n}_{k}}), D(\widehat{\beta^{n}_{k}}))=0.$

Now we will prove that top degree monomials of $P$ are zero :
$$ P(Y, Z)=\sum^{d} \limits_{i=0}A_{d-i, i}(s, s^{-1})Y^{d-i}Z^{i}+\sum \limits_{0\leq i+j\leq d-1}A_{i,j}(s, s^{-1})Y^{i}Z^{j} $$
For $Y=V(\widehat{\beta^{n}_{k}})$ and $Z=D(\widehat{\beta^{n}_{k}})$ we have
 $$\Big[\sum^{d} \limits_{i=0}\frac{1}{2^{(n-k)i}}A_{d-i,i}(s, s^{-1})\Big]\Big(\frac{s^2+1}{s}\Big)^{d(n-k)}+ \ldots =0$$ and for $n-k$ big enough we have $$\max (\deg A_{d-i, i})+d(n-k)> \max (\deg A_{i, j})+(d-1)(n-k),$$
therefore Laurent polynomial $\sum^{d} \limits_{i=0}\frac{1}{2^{(n-k)i}}A_{d-i, i}(s, s^{-1})$ should be zero.
Induction on $i$ ends the proof : for $n-k \rightarrow \infty$ we obtain $A_{d, 0}(s, s^{-1})=0,$ and if $A_{d,0}=A_{d-1, 1}= \ldots = A_{d-i+1, i-1}=0$, then $2^{i(n-k)}[\frac{1}{2^{i(n-k)}}A_{d-i, i}(s, s^{-1})+\frac{1}{2^{(i+1)(n-k)}}A_{d-i-1, i+1}(s, s^{-1})+ \ldots]$ is zero and $n-k\rightarrow \infty$ gives $A_{d-i, i}(s, s^{-1})=0$.
\end{proof}
As a final example we take $\mathcal{A}=\{\hbox{isotopy classes of links}\}$, $f_{1}(L)=\nabla(L)$, $f_{2}(L)=V(L)$ and $f_{3}(L)=D(L)$ and we show that $\{\nabla(L),V(L),D(L)\}_{\mathcal{A}}$ is algebraically independent over $\mathbb{C}[s,s^{-1}]$.\\
\emph{Proof of Theorem \ref{algebraically independent}} Take $p$ a prime number, $n$ an even natural number and the link $L_{p,n}=\widehat{\beta}_{p,n}$, where $\beta_{p,n}$ is the $n+1$-braid $x_{1}^{2p+1}x_{2}^{2p+1}\ldots x_{n}^{2p+1}$. The corresponding polynomials of $L_{p,n}$ have leading terms $\nabla(L_{p,n})=(s^{2p}+\ldots )^{n}$, $V(L_{p,n})=(s^{6p+2}+\ldots)^{n}$, and $D(L_{p,n})=(ps^{2p+2}+\ldots)^{n}$. Suppose  that
 $$P(X,Y,Z)=\sum \limits_{i,j,k}A_{i,j,k}(s,s^{-1})X^{i}Y^{j}Z^{k}$$
 is a nonzero polynomial with coefficients in $\mathbb{C}[s,s^{-1}]$ which is zero for $(X,Y,Z)=(\nabla(L_{p,n}),V(L_{p,n}),D(L_{p,n}))$ (only nonzero coefficient $A_{i,j,k}$ are in this sum, $\alpha_{i,j,k}s^{a_{i,j,k}}$is the leading coefficient  of $A_{i,j,k}(s,s^{-1})$). We will find a contradiction if $P$ is nonzero. The leading term of the polynomial in $s$ $P(L_{p,n})=P(\nabla(L_{p,n}),V(L_{p,n}),D(L_{p,n}))$ is given by
$$\sum\alpha_{i,j,k}p^{nk}s^{a_{i,j,k}+2pni+(6p+2)nj+(2p+2)nk}$$ where $a_{i,j,k}+2pni+(6p+2)nj+(2p+2)nk$ is maximal ($i,j,k$ are variable). Its coefficient $\sum\limits_{i,j,k}\alpha_{i,j,k}p^{nk}$ should be zero, and the polynomial equation (complex coefficients) $\sum\limits_{i,j,k}\alpha_{i,j,k}t^{k}=0$ cannot have infinitely many solutions $(p^{n})$, therefore for a fixed $k$, $\sum\limits_{i,j}\alpha_{i,j,k}=0$. From maximality condition $i$ and $j$ are such that the sum $S_{i,j}=a_{i,j,k}+2pni+(bp+2)nj$ is constant in $i$ and $j$ (for any  $p$ and $n$); we show that this implies the uniqueness of the pair $(i,j)$, hence we find a contradiction: $\alpha_{i,j,k}=0$. First choose $n$ greater than any difference $|a_{i,j,k}-a_{h,l,k}|$, next choose $p$ greater than any difference $|a_{i,j,k}-a_{h,l,k}|$ and  $|a_{i,j,k}+2nj-a_{h,l,k}-2nl|$ (condition on $n$ shows that the last difference is not $0$ for $j\neq l$). If there are two pairs $(i,j),(h,j)$ with $S_{i,j}=S_{h,j}$, we find $a_{i,j,k}+2pni=a_{h,j,k}+2pnh$, hence $p\,|\, (a_{i,j,k}-a_{h,j,k})$ and this is possible only for $a_{i,j,k}=a_{h,j,k}$, and the previous equation gives $i=h$. If there are two pairs $(i,j),(h,l)$ ($j\neq l$) with $S_{i,j}=S_{h,l}$, we find $a_{i,j,k}+2pni+(6p+2)nj=a_{h,l,k}+2pnh+(6p+2)nl$, hence $p\,|\,(a_{i,j,k}+2nj-a_{h,l,k}-2nl)$ and this is not possible because $0<|a_{i,j,k}+2nj-a_{h,l,k}-2nl|<p$.
 \begin{flushright}
$\Box$
\end{flushright}
\begin{example}(Eliahou, Kauffman, Thislethwaite \cite{trivial: Jones}) The link $LL_{2}(2)$ has the same Jones and Alexander-Conway  polynomials, and different $D$ polynomials. Using the relative expansion formula $\mathcal{L}(-5,6)$ we have to compute $D_{\mathcal{L}}(0,0)$,  $D_{\mathcal{L}}(0,1)$,  $D_{\mathcal{L}}(1,0)$ and $D_{\mathcal{L}}(1,1)$:
 $$D_{\mathcal{L}}(-5,6)=-30sD_{\mathcal{L}}(0,0)+36D_{\mathcal{L}}(0,1)+25D_{\mathcal{L}}(1,0)-30s^{-1}D_{\mathcal{L}}(1,1).$$
 $\mathcal{L}(1,0)$ is the unlink $\bigcirc\bigcirc\bigcirc$,  $\mathcal{L}(0,0)$ is Hopf$^{-}\coprod\bigcirc$ and  $\mathcal{L}(1,1)$ is Hopf$^{+}\coprod\bigcirc$ where Hopf$^{-}=\widehat{x_{1}^{2}}$ and Hopf$^{+}=\widehat{x_{1}^{-2}}$. $\mathcal{L}(0,1)$ is a four component link: after two skein decompositions we find
Hopf$^{+}\coprod\bigcirc$, the unlink $\bigcirc\bigcirc\bigcirc$ and closure of the $4$-braid $x_{1}^{-1}x_{2}^{2}x_{1}^{-1}x_{3}^{-3}x_{2}^{2}x_{3}$. Computing only the leading terms, the main contribution comes from the first link and we have
$D_{\mathcal{L}}(-5,6)=\frac{15}{2}s^{5}+\ldots$ and $D(\bigcirc\bigcirc)=\frac{s^{2}+1}{2s}$.
\end{example}

\textbf{Acknowledgements:} The first author has presented a preliminary version of this paper at the  Knot Theory Conference held at ICTP in May, 2009.  She is grateful to ICTP for inviting her and giving the opportunity to present these results and also to Professor Hugh Morton and Professor Jozef Przytycki for their valuable suggestions they gave during that conference. Both  authors are thankful to Pakistan Higher Education Commission for supporting  this research.

\medskip
%\textbf{2000 Mathematics Subject Classification: }%05A15, 11B39

\end{document}